\documentstyle[12pt,amssymbols]{article}

\newfont{\g}{eufm10 scaled\magstep1}
\newfont{\gs}{eufm7}

\newcommand{\Sym}{\mbox{\g S}}

\newfont{\bg}{cmr10 scaled\magstep4}

\newcommand{\bigzerou}{\smash{\hbox{\lower1.7ex\hbox{\bg 0}}}}

\def\N{\Bbb N}
\def\C{\Bbb C}

\def\Q{\Bbb Q}

\def\I{{\cal I}}

\def\SP{{\cal S\!P}}
\def\P{{\cal P}}
\def\<{\langle}
\def\>{\rangle}

\def\a{\alpha}

\def\ga{\gamma}
\def\Ga{\Gamma}
\def\e{\epsilon}

\def\l{\lambda}
\def\L{\Lambda}
\def\d{\delta}
\def\le{\lambda^0}
\def\lo{\lambda^1}

\def\pr{\prod}
\def\ol{\overline}

\def\rbar{\bar{r}}

\def\sgn{\mbox{\rm sgn}}

\newenvironment{proof}%
{\noindent{\bf Proof. }}%
{\hspace{12pt}$\blacksquare$\par\vspace{12pt}}

\newcounter{sect}
\def\thesect{\arabic{sect}}
\newtheorem{definition}{Definition}[sect]
\newtheorem{theorem}[definition]{Theorem}
\newtheorem{lemma}[definition]{Lemma}
\newtheorem{proposition}[definition]{Proposition}

\begin{document}
\setcounter{page}{0}
\thispagestyle{empty}
\vspace{1cm}
\begin{large}
\begin{center}
{\bf ON REDUCED $Q$-FUNCTIONS}
\end{center}
\end{large}

\vspace{1.0cm}

\begin{center}
Tatsuhiro NAKAJIMA

{\it Department of Physics, Tokyo Metropolitan University}\\
{\it Minami-Ohsawa 1-1, Hachioji-shi, Tokyo 192-03, Japan}\\
and

Hiro-Fumi YAMADA

{\it Department of Mathematics, Tokyo Metropolitan University}\\
{\it Minami-Ohsawa 1-1, Hachioji-shi, Tokyo 192-03, Japan}
\vspace{1.0cm}

{\it To Kiyosato Okamoto with affection and admiration}

\end{center}
\vspace{2cm}
\begin{center}
{\bf Abstract}
\end{center}
Schur's $Q$-functions with reduced variables are discussed by
employing a combinatorics of strict partitions.
They are called reduced $Q$-functions.
We give a description of the linear relations among reduced $Q$-functions.

\newpage
\setcounter{sect}{0}
\begin{flushleft}
{\large\bf \S \thesect \hspace{5mm} Introduction}
\end{flushleft}
\vspace{3mm}

$Q$-functions were introduced by Schur in his study of projective
representations of symmetric groups.
They are symmetric functions and, if we express them in terms of
the power sum symmetric functions, each coefficient essentially gives
the spin character of the symmetric group.

This note deals with the ``$r$-reduced'' $Q$-functions which are
defined by putting $p_{jr}=0$ for $j=1$, $3$, $5,\ldots$ in the power
sum expression.
When $r=p$ is a prime number, they play a role in $p$-modular projective
representations of symmetric groups.

In a previous work [NY] we showed that $r$-reduced $Q$-functions
are weight vectors of the basic representation of the affine Lie algebra
$A^{(2)}_{2t}$ $(r=2t+1)$ and chose a proper basis for each weight space.

To be more precise, let $\a_i$ (resp. $\a^{\vee}_i$) $(0\leq i\leq t)$
be the simple roots (resp. coroots) of the affine Lie algebra $A^{(2)}_{2t}$
and $\d=2\sum_{i=0}^{t-1} \a_i+\a_t$ be the fundamental imaginary root.
The irreducible representation with highest weight $\L_0$ is called the
basic representation, where $\L_0(\a_i^{\vee})=\d_{i0}$.
The set of weights is described by
$$
P=\{W\L_0-n\d\ ;\ n\in \N \},
$$
where $W$ is the Weyl group.
This basic representation can be realized on the polynomial ring
$\C[t_j\ ;\ j\geq 1,\ \mbox{odd and }j\not\equiv 0\ (mod\ r)]$.
In this realization each $r$-reduced $Q$-function turns out
to be a weight vector.
We answered in [NY] the question to which weight the given $r$-reduced
$Q$-function belongs, by using a combinatorics of strict partitions,
i.e., bar cores and bar quotients.

By virtue of the above result we further investigate the $r$-reduced
$Q$-functions themselves.
The result in this note is an explicit description of the linear
relations satisfied by $r$-reduced $Q$-functions.
We remark that such a description has been obtained for ``$r$-reduced
Schur functions'' [ANY1,2].

\stepcounter{sect}
\begin{flushleft}
{\large\bf \S \thesect \hspace{5mm} Combinatorics of strict partitions}
\end{flushleft}
\vspace{3mm}

We first present a collection of definitions and results
concerning strict partitions and Schur's $Q$-functions,
mainly referring to the new book of Macdonald [M].

The set of all partitions of $n$ is denoted by $\P_n$ and
the set of all partitions by $\P$.
A partition $\l=(\l_1,\ldots,\l_l)$ is said to be strict if
$\l_1>\l_2>\cdots>\l_l>0$.
The set of all strict partitions of $n$ is denoted by $\SP_n$ and
the set of all strict partitions by $\SP$.

To each element $\l\in \SP$ we can associate the shifted diagram
$S(\l)$, which is obtained from the usual Young diagram of $\l$
by shifting the $i$-th row $(i-1)$-positions to the right.
The $j$-th square in the $i$-th row will be called the $(i,j)$-square.
For $i$, $j$ and $k$ with $k\geq i$, we define an $(i,j)_k$-bar of $\l$
as the set of squares $(a,b)\in S(\l)$ satisfying
(1) if $k>i$, then $a=i$ or $k$,
(2) if $k=i$, then $a=i$, $b\geq j$ and $j-1\neq \l_{a'}$ for $a'\geq i+1$.
The length of a bar is a number of squares it contains.
A bar of length $r$ will be often called an $r$-bar.
Throughout this note we always assume that $r$ is an odd number.

A strict partition is an $r$-bar core if it contains no $r$-bars.
Given a strict partition $\l$ we obtain another strict partition
by removing an $r$-bar and rearranging the rows in descending order.
By repeating this process as often as possible we shall end up with an
$r$-bar core $\l^c$, which is called the $r$-bar core of $\l$.

For $0\leq i \leq r-1$ we define sequences
$$
X^i(\l)=\{x\in \N ; \l_k=rx+i \mbox{ for some } k\}.
$$
The $r$-bar quotient $\l^q=(\l^0,\l^1,\ldots,\l^t)$ $(t=(r-1)/2)$is a
$(t+1)$-tuple of partitions given by
$\l^0=X^0(\l)$ and $\l^i=(X^i(\l)|X^{r-i}(\l))$ $(1\leq i\leq t)$
in Frobenius notation.
Note that $\l^0$ is a strict partition.
We remark that any $r$-bar core and $r$-bar quotient
uniquely determines a strict partition [O].
In view of this one-to-one correspondence we often write $(\l^c,\l^q)$
instead of $\l$, e.g., $Q_{(\l^c,\l^q)}(x)$, etc.

We introduce $r$-bar sign for strict partitions.
Let $\bar{H}$ be the $(i,j)_k$-bar of length $r$ in $\l\in \SP$.
The leglength of $\bar{H}$ is defined by
$$
b(\bar{H})=\left\{
\begin{array}{ll}
\l_k+\#\{a ; \l_i>\l_a>\l_k\} & \mbox{if\hspace{12pt}} k>i \\
\#\{a ; \l_i>\l_a>\l_i-r\} & \mbox{if\hspace{12pt}} k=i.
\end{array}
\right.
$$
If the strict partition $\mu$ is obtained from $\l$ by removing
$r$-bars $H_1$, $H_2,\ldots, H_q$ successively, then the $r$-bar sign
of $\l$ relative to $\mu$ is defined by
$$
\d_{\rbar}(\l,\mu)=\pr_{m=1}^q (-1)^{b(\bar{H}_m)}.
$$
We set $\d_{\rbar}(\l)=\d_{\rbar}(\l,\l^c)$ and call it the $r$-bar sign
of $\l$.
We should remark that $\d_{\rbar}(\l,\mu)$ does not depend on the choice of
$r$-bars being removed in going from $\l$ to $\mu$ [MO].

Let $\L_{\Q}$ denote the graded $\Q$-algebra of symmetric functions in
countably many independent variables $x=(x_1$, $x_2, \ldots)$.
Among various bases of $\L_{\Q}$, we just mention here that the
Schur functions $s_{\l}(x)$ $(\l\in P)$ are convenient for our point of view.
Let $\Ga_{\Q}$ be the subalgebra of $\L_{\Q}$ generated by
$1$, $p_1$, $p_3$, $p_5, \ldots$;
$$
\Ga_{\Q}=\Q[p_1,p_3,p_5,\ldots],
$$
where $p_j=x_1^j+x_2^j+\cdots$ is the $j$-th power sum.

Schur's $Q$-functions $Q_{\l}(x)$ $(\l\in \P)$ are defined by
specializing to $t=-1$ in the Hall-Littlewood symmetric functions
$Q_{\l}(x; t)$.
This specialization implies that $Q_{\l}(x)=0$ unless $\l\in \SP$.
It is known that $\{Q_{\l}(x) ; \l\in \SP\}$ gives a basis of $\Ga_{\Q}$.
When $\mu =\sigma\l=(\l_{\sigma(1)},\ldots,\l_{\sigma(l)})$
$(\sigma\in \Sym_l)$ is a derangement of a strict partition $\l$,
then we define $Q_{\mu}(x)=\sgn(\mu) Q_{\l}(x)$, where we have put
$\sgn(\mu)=\sgn(\sigma)$.

The $Q$-functions are related to the power sums via the character
values of the symmetric group.
An explicit expression is given as follows.
Let $\zeta_{\l}(\pi)$ ($\l\in \SP_n$) be the value of the irreducible
projective character $\zeta_{\l}$ of the symmetric group $\Sym_n$
at the class $\pi=(1^{\pi_1} 3^{\pi_3} 5^{\pi_5} \cdots)$ and
set $z_{\pi}=\prod_{j\geq 1} \pi_j! j^{\pi_j}$.
Then we have
$$
Q_{\l}(x)=\sum_{\pi} 2^{[(\ell(\l)+\ell(\pi)+1)/2]} z_{\pi}^{-1}
          \zeta_{\l}(\pi) p_1^{\pi_1} p_3^{\pi_3} p_5^{\pi_5} \cdots,
$$
where the summation runs over all odd partitions
$\pi=(1^{\pi_1} 3^{\pi_3} 5^{\pi_5} \cdots)$ of size $n$, and
$\ell(\l)$ denotes the length of the partition $\l$.

The following formula is of importance for our purpose.
\begin{theorem}
Let $j\geq 1$ be odd and $\l=(\l_1,\ldots,\l_l)\in \SP$. Then
$$
p_j Q_{\l}(x)=\frac{1}{2}\sum_{i=1}^{(j-1)/2} (-1)^i Q_{(\l;j-i,i)}(x)
             +\frac{1}{2} Q_{(\l;j)}(x)
             +\sum_{i=1}^l Q_{\l+j\e_i}(x),
$$
where $(\l;j-i,i)=(\l_1,\ldots,\l_l,j-i,i)$, $(\l;j)=(\l_1,\ldots,\l_l,j)$
and $\e_i=(\d_{ik})_{1\leq k\leq l}$.
\end{theorem}
This is a direct consequence of [J, Theorem 2.4] by looking at
the correspondence between removal of {\it strips/double strips} and
that of $j$-bars.

\stepcounter{sect}
\begin{flushleft}
{\large\bf \S \thesect \hspace{5mm}
Reduced $Q$-functions and their linear relations}
\end{flushleft}
\vspace{3mm}

Let
$\Ga^{(r)}_{\Q}
=\Q[\,p_j \,; j\geq 1, \mbox{ odd and } j\not\equiv 0\mbox{ (mod $r$)}]$
and define the $r$-reduced $Q$-functions by
$$
Q^{(r)}_{\l}(x)
=\left.Q_{\l}(x)\right|_{p_r=p_{3r}=p_{5r}=\cdots=0}\in \Ga^{(r)}_{\Q}.
$$

\begin{proposition}
Let $\l$ be a strict partition and  $\l^q=(\le,\ldots,\l^t)$ $(t=(r-1)/2)$
be the $r$-bar quotient of $\l$. Then
the set $\{Q^{(r)}_\l(x) ; \le=\emptyset\}$ gives a basis of $\Ga^{(r)}_{\Q}$.
\end{proposition}
To prove this proposition we need the following lemma.
\begin{lemma}
Let $V$ be an infinite dimensional $\Q$-vector space with basis
$\{v_{\l} ; \l\in \SP\}$ and define $v_{\sigma\l}=\sgn(\sigma) v_{\l}$ for
permutations $\sigma$; $v_{\l}=0$ if $\l\not\in \SP$.
\begin{enumerate}
\item If we set
$$
V_1:=V\left/\sum_{\l\in \SP} \Q \left(\sum_{j\geq 1, {\rm odd}} \left(
\frac{1}{2} \sum_{i=1}^{(j-1)/2} (-1)^i v_{(\l;j-i,i)}
+\frac{1}{2} v_{(\l;j)}
+\sum_{i=1}^{\ell(\l)} v_{\l+j\e_i} \right.
\right) \right),
$$
then $V_1 \cong \Q$.
\item If we set
$$
V_r:=V\left/\sum_{\l\in \SP} \Q \left(\sum_{j\geq 1,{\rm odd}} \left(
\frac{1}{2} \sum_{i=1}^{(rj-1)/2} (-1)^i v_{(\l;rj-i,i)}
+\frac{1}{2} v_{(\l;rj)}
+\sum_{i=1}^{\ell(\l)} v_{\l+rj\e_i} \right.
\right) \right),
$$
then $V_r\cong \Ga^{(r)}_{\Q}$.
\end{enumerate}
\end{lemma}

\begin{proof}
(1) There is a canonical linear surjection $\ga$ from $\Ga_{\Q}$ to
$V_1$ which maps $Q_{\l}(x)$ to $v_{\l}$ for $\l\in \SP$.
The kernel of $\ga$ coincides with the maximal ideal
$I=(p_1,p_3,p_5,\ldots)$ because of Theorem 1.1.
Since the algebra $\Ga_{\Q}/\I$ is isomorphic to $\Q$, so is $V_1$.

\noindent
(2) Consider a linear surjection
$\ga^{(r)}: \Ga_{\Q} \longrightarrow V_r$ defined by
$\ga^{(r)}(Q_\l(x))=v_{\l}$ for $\l\in \SP$.
If we set $\I^{(r)}=(p_r,p_{2r},p_{3r},\ldots)$, then, by
Theorem 1.1, we see that $\ga^{(r)}(\I^{(r)})=0$ and can define a
linear surjection $\bar{\ga}^{(r)}: \Ga_{\Q}/\I^{(r)} \longrightarrow V_r$.
On the other hand we can define a linear surjection
$\eta: V_r \longrightarrow \Ga_{\Q}^{(r)}$ by $\eta(v_{\l})=Q^{(r)}_{\l}(x)$
for $\l\in \SP$. The composition
$\eta \circ \bar{\ga}^{(r)}: \Ga_{\Q}/\I^{(r)} \longrightarrow \Ga_{\Q}^{(r)}$
gives a linear isomorphism.
Hence $\eta$ is a linear isomorphism as desired.
\end{proof}

\noindent
{\bf Proof of Proposition 2.1. } We give a proof only for the case
$r=3$ in order to avoid the too complicated notation.
The general case can be shown in a same fashion.

First we shall see that the set $\{Q^{(3)}_{\l} ; \le=\emptyset\}$
spans $\Ga^{(3)}_{\Q}$.
We introduce a filtration on $\Ga^{(3)}_{\Q}$
by $F_n\Ga^{(3)}_{\Q}=\displaystyle{\sum_{\l; |\le|\leq n}}
\Q Q^{(3)}_{(\l^c,\le,\lo)}(x)$
and the associated graded module
$\ol{F}(\Ga^{(3)}_{\Q}):=
\displaystyle{\bigoplus_{n\geq 0} F_n\Ga^{(3)}_{\Q}/F_{n-1}\Ga^{(3)}_{\Q}}$,
where $F_{-1}\Ga^{(3)}_{\Q}=\{0\}$.
For any positive odd integer $j$ the $3$-reduced $Q$-functions
$Q^{(3)}_{\l}(x)$ satisfy
\begin{eqnarray*}
0&=&\frac{1}{2} \sum_{i=1}^{(3j-1)/2} Q^{(3)}_{(\l; 3j-i,i)}(x)
+\frac{1}{2}  Q^{(3)}_{(\l; 3j)}(x)
+\sum_{i=1}^{\ell(\l)} Q^{(3)}_{\l+3j\e_i}(x) \\
&=&\frac{1}{2} \sum_{i=1}^{(j-1)/2}  Q^{(3)}_{(\l^c,(\le;j-i,i), \lo)}(x)
+\frac{1}{2}  Q^{(3)}_{(\l^c,(\le;j) ,\lo)}(x)
+\sum_{i=1}^{\ell(\le)}  Q^{(3)}_{(\l^c,\le+j\e_i, \lo)}(x) \\
&+&\frac{1}{2} \sum_{i\not\equiv 0 \,({\rm mod}\, 3)}
Q^{(3)}_{(\l^c,\le, (\l;3j-i,i)^1)}
+\sum_{i=1}^{\ell(\lo)} Q^{(3)}_{(\l^c,\le,\lo+j\e_i)}(x).
\end{eqnarray*}
If we write $\ol{Q}^{(3)}_{\l}(x)=Q^{(3)}_{\l}(x)+F_{n-1}\Ga^{(3)}_{\Q}$ for
$|\le|=n$, then we see that in $\ol{F}(\Ga^{(3)}_{\Q})$
$$
\sum_{i=1}^{\ell(\le)}  \ol{Q}^{(3)}_{(\l^c,\le+j\e_i, \lo)}(x)
+\frac{1}{2}  \ol{Q}^{(3)}_{(\l^c, (\le;j) ,\lo)}(x)
+\frac{1}{2} \sum_{i=1}^{(j-1)/2}  \ol{Q}^{(3)}_{(\l^c, (\le;j-i,i), \lo)}(x)
=0.
$$
By applying Lemma 2.2 (1)  it turns out that
$\displaystyle{\sum_{\le\in \SP} \Q \ol{Q}^{(3)}_{(\l^c,\le,\lo)}(x)}
=\Q \ol{Q}^{(3)}_{(\l^c,\emptyset,\lo)}(x)$.
This proves that $\{Q^{(3)}_{\l}(x) ; \le=\emptyset\}$ spans $\Ga^{(3)}_{\Q}$.

By looking at the one-to-one correspondence between the strict partitions
and the tuples of $r$-cores and $r$-quotients, we see that the number of
those strict partitions $\l\in \SP_n$ such that $\le=\emptyset$ coincides
with the dimension of the subspace of homogeneous polynomials of
degree $n$ in $\Ga_{\Q}^{(r)}$, where we count $\deg p_j=j$.
This completes the proof.
\hspace{12pt}$\blacksquare$
\vspace{12pt}

We shall state the linear relations satisfied by the $r$-reduced
$Q$-functions.
To this end, we define two types of Littlewood-Richardson-like coefficients
$f_{\l_1\cdots \l_t}^{\nu}$ $(\l_1,\ldots,\l_t,\nu\in \SP)$
and $h_{\l \mu}^{\nu}$ $(\l\in \SP,\, \mu, \nu\in \P)$, respectively by
\begin{eqnarray*}
& &P_{\l_1}(x)\cdots P_{\l_t}(x)
=\sum_{\nu\in \SP} f_{\l_1\cdots \l_t}^{\nu} P_{\nu}(x), \\
& &P_{\l}(x) s_{\mu}(x)=\sum_{\nu\in \P} h_{\l \mu}^{\nu }s_{\nu}(x),
\end{eqnarray*}
where we set $P_{\l}(x)=2^{-\ell(\l)} Q_{\l}(x)$ for $\l\in \SP$.
Note that $f^{\nu}_{\l_1\cdots \l_t}$ and $h_{\l \mu}^{\nu}$
are non-negative integers.

\begin{theorem}
Let $\l$ be a strict partition, $\l^c$ and $\l^q=(\le,\lo,\ldots,\l^t)$ be its
$r$-bar core and $r$-bar quotient, respectively. Then
$$
Q^{(r)}_{\l}(x)=2^{[(\ell(\l)+\ell(\le))/2]} (-1)^{|\le|}
\d_{\bar{r}}(\l) \sum_{\mu, \nu_1,\ldots,\nu_t}
2^{-[\ell(\mu)/2]} \d_{\bar{r}}(\mu)
f^{\le}_{\nu_1\cdots \nu_t}
h_{\nu_1 \lo}^{\mu^1}\cdots h_{\nu_t \l^t}^{\mu^t} Q_{\mu}^{(r)}(x),
$$
where summation runs over the strict partitions $\nu_1, \ldots,\nu_t$ and
$\mu$ such that $|\mu|=|\l|$, $\mu^c=\l^c$ and $\mu^0=\emptyset$.
\end{theorem}

\begin{proof}
Again we give a proof only for the case $r=3$.
For the general case we should notice the following identity :
$$
Q_{\nu}(x^{(1)},x^{(2)},\ldots,x^{(t)})=
\sum_{\l_1,\ldots,\l_t} f_{\l_1\cdots \l_t}^{\nu}
Q_{\l_1}(x^{(1)}) Q_{\l_2}(x^{(2)})\cdots Q_{\l_t}(x^{(t)}),
$$
where summation runs over the strict partitions $\l_1,\ldots,\l_t$.
Here $(x^{(1)},x^{(2)},\ldots,x^{(t)})$ are $t$ sets of variables and
$Q_{\nu}(x^{(1)},\ldots,x^{(t)})$ denotes the $Q$-function in the
set of variables \par\noindent
$(x^{(1)}_1,x^{(1)}_2,\ldots,x^{(2)}_1,x^{(2)}_2,\ldots,
x^{(t)}_1,x^{(t)}_2,\ldots)$.

Set
$$
F_{\l}(x)=\sgn(\l) 2^{[(\ell(\l)+\ell(\le))/2]} (-1)^{|\le|}
\d_{\bar{3}}(\l) \sum_{\mu} 2^{-[\ell(\mu)/2]} \d_{\bar{3}}(\mu)
h_{\le \lo}^{\mu^1}Q_{\mu}^{(3)}(x).
$$
We shall prove that
$$
\sum_{i=1}^{(3j-1)/2}(-1)^i F_{(\l;3j-i,i)}(x)
+\frac{1}{2} F_{(\l;3j)}(x)+\sum_{i=1}^{\ell(\l)} F_{\l+3j\e_i}(x)=0
$$
for any positive odd integer $j$.
For our purpose it is enough to show that the coefficient of
$Q_{\mu}^{(3)}(x)$ vanishes for each $\mu\in \SP$.
Since corollary in [MO, (3.8)] leads to
\begin{eqnarray*}
\lefteqn{\d_{\bar{3}}(\tilde{\l},\l)} \\
&=& \left\{
\begin{array}{ll}
\sgn(\tilde{\l})\sgn(\tilde{\l}^k)&
{\rm for}\ \tilde{\l}=\l+3j\e_i\ (\,k=0,1,\ i=1,\ldots,\ell(\l))  \\
\sgn(\tilde{\l})\sgn(\tilde{\l}^0) &
{\rm for}\ \tilde{\l}=(\l; 3j)\ {\rm and}\ (\l; 3j-3j',3j')\
\left(\, j'\geq \frac{j-1}{2}\right) \\
(-1)^i \sgn(\tilde{\l}) \sgn(\tilde{\l}^1)  &
{\rm for}\ \tilde{\l}=(\l; 3j-i,i)\
\left(\,i\geq\frac{3j-1}{2}, i\not\equiv 0 ({\rm mod}\ r)\right),
\end{array}
\right.
\end{eqnarray*}
our equation reads
\begin{eqnarray*}
& &\frac{1}{2}\sum_{i< j/2} 2^2(-1)^i \sgn((\le; j-i,i)) (-1)^j
h_{(\le; j-i,i) \lo}^{\mu^1} \\
&+&\frac{1}{2}\!\!\sum_{i\not\equiv 0 \,\,({\rm mod}\,\,3)}\!\!\!
2\, \sgn((\l; 3j-i,i)^1)
h_{\le (\l; 3j-i,i)^1}^{\mu^1}
+\frac{1}{2} \sgn((\le; j)) 2(-1)^j h_{(\le; j) \lo}^{\mu^1} \\
&+&\sum_{i} \sgn((\l+3j\e_i)^0) (-1)^j
h_{(\l+3j\e_i)^0 \lo}^{\mu^1}
+\sum_{i} \sgn((\l+3j\e_i)^1) h_{\le (\l+3j\e_i)^1}^{\mu^1} \\
&=& 0.
\end{eqnarray*}
Multiplying both sides by $s_{\mu^1}(x)$ and taking a summation over $\mu^1$,
we see that the equation reduces to
\begin{eqnarray*}
& &(-1)^j \left\{\frac{1}{2} \sum_{i<j/2} (-1)^i Q_{(\le; j-i,i)}(x) s_{\lo}(x)
+\frac{1}{2}Q_{(\le; j)}(x) s_{\lo}(x)+
\sum_{i} Q_{\le+j\e_i}(x) s_{\lo}(x)\right\}  \\
&+&\left\{ \sum_{i} Q_{\le}(x) s_{\lo+j\e_i}(x)+\!\!\!\!\!
\sum_{\scriptstyle i<j/2\atop\scriptstyle i\not\equiv 0\,({\rm mod}\,3)}
\!\!\!\!\!\!
Q_{\le}(x) s_{(\l,3j-i,i)^1}(x)\right\}=0.
\end{eqnarray*}
Note that there is a one-to-one correspondence between the $3j$-bars
in $\l$ and the $j$-bars in $\l^q=(\le,\lo)$ (a $j$-bar in $\l^q$ means
a $j$-bar in $\le$ or a $j$-hook in $\lo$). Hence we have
$$
-(p_j Q_{\le}(x))s_{\lo}(x)+Q_{\le}(x)(p_j s_{\lo}(x))=0.
$$

By Lemma 2.2 (3) there exists a linear isomorphism sending $F_{\l}(x)$
to $Q^{(3)}_{\l}(x)$ and it turns out to be the identity since
$Q^{(3)}_{\l}(x)$ equals $F_{\l}(x)$ for the basis given in Proposition 2.1.
\end{proof}

\end{document}